\documentclass[12pt, reqno]{amsart}
\usepackage{amsmath, amsfonts, amssymb}
\usepackage{hyperref}
\usepackage[cp1251]{inputenc}
\hypersetup{colorlinks=true,linkcolor=red, anchorcolor=green, citecolor=cyan, urlcolor=red, filecolor=magenta, pdftoolbar=true}

\newtheorem{theorem}{Theorem}[section]

\theoremstyle{definition}
\newtheorem{definition}{Definition}[section]

\theoremstyle{remark}
\newtheorem{remark}{Remark}[section]
\numberwithin{equation}{section}

\begin{document}

\setcounter{page}{1}

\title[$p$-Laplace heat equations with combined nonlinearities]{Critical exponents for the $p$-Laplace heat equations with combined nonlinearities}

\author[B. T. Torebek]{Berikbol T. Torebek}

\address{\textcolor[rgb]{0.00,0.00,0.84}{Berikbol T. Torebek \newline Department of Mathematics: Analysis, Logic and Discrete Mathematics \newline Ghent University, Krijgslaan 281, B9000 Ghent, Belgium
\newline and \newline Institute of Mathematics and Mathematical Modeling \newline 125 Pushkin str., 050010 Almaty, Kazakhstan}}
\email{\textcolor[rgb]{0.00,0.00,0.84}{berikbol.torebek@ugent.be }}




\thanks{This research has been/was/is funded by the Science Committee of the Ministry of Science and Higher Education of the Republic of Kazakhstan (Grant No. AP14869090) and by the FWO Odysseus 1 grant G.0H94.18N: Analysis and Partial Differential Equations and by the Methusalem programme of the Ghent University Special Research Fund (BOF) (Grant number 01M01021). No new data was collected or generated during the course of research.}

\begin{abstract}
This paper studies the large-time behavior of solutions to the quasilinear inhomogeneous parabolic equation with combined nonlinearities. This equation is a natural extension of the heat equations with combined nonlinearities considered by Jleli-Samet-Souplet \cite{Jleli1} (Proc AMS, 2020).
Firstly, we focus on an interesting phenomenon of discontinuity of the critical exponents. In particular, we will fill the gap in the results of \cite{Jleli1} for the critical case. We are also interested in the influence of the forcing term on the critical behavior of the considered problem, so we will define another critical exponent depending on the forcing term.
\end{abstract} \maketitle
\tableofcontents
\section{Introduction}
This paper is concerned with the large time behavior of global weak solutions to the following degenerate heat equation with a forcing term
\begin{equation}\label{1}
\left\{\begin{array}{l}u_t-\Delta_p u=\lambda|u|^\alpha+\mu |\nabla u|^\beta+f(x),\, (t,x)\in (0,T)\times\mathbb{R}^n, T>0,\\{}\\u(0,x)=u_0(x),\, x\in \mathbb{R}^n,\end{array}\right.
\end{equation} where $\lambda, \mu\in\mathbb{R},$ $\alpha,\beta>0,$ $\nabla$ is a gradient operator and $$\Delta_p u=\text{div}\left(|\nabla u|^{p-2}\nabla u\right)\,\left(p> \frac{2n}{n+1}\right)$$ is the $p$-Laplace operator. Specifically, we are interested in global and blow-up solutions to the considered problem.

If $p=2$ and $\lambda = 0,$ then the equation \eqref{1} is called a viscous Hamilton-Jacobi equation, which has the form
$$u_t-\Delta u=\mu |\nabla u|^\beta+f(x),$$
and it appears in the physical theory of growth and roughening of surfaces, where it is known as the
deterministic version of the KPZ (Kardar-Parisi-Zhang) model (when $\lambda = 0,$ $p=2$ and $\beta=2$), describing the profile of a growing interface in certain physical models (see \cite{Kardar, Krug}).

The local existence of problems of type \eqref{1} has been studied by many authors (for example, see \cite{Serrin, Vespri, Lad, Lian, Lu, Shang} and the references therein). For example, the local existence of more general cases of problem \eqref{1} was studied in \cite{Serrin, Vespri, Lad}, and in \cite{Lian, Lu, Shang} the existence of the local solution of problem \eqref{1} was proved for locally integrable initial data.

The main purpose of this paper is to study the critical exponents of the
existence and nonexistence of global solutions to problem \eqref{1}.

Namely, we are interested in the following questions:
\begin{itemize}
  \item Do the critical exponents of problem \eqref{1} remain discontinuous for all $p>\frac{2n}{n+1}$? It was demonstrated in \cite{Jleli1} that when $p=2$, the critical exponents are discontinuous.
  \item  Does the considered problem admit global or blow-up solutions in the critical cases of $\alpha$ and $\beta$? This question is still open for any $p>\frac{2n}{n+1}$.
\end{itemize}

More detailed information about the history of critical exponents for evolution equations related to problem \eqref{1} is given in Subsection \ref{Hist}.

\subsection{Main results}
We will give our first result by the following theorem. Here and below, we use the convention $\alpha_{cr} = \infty$ for $n \leq p.$

We first give the definition of weak solutions to problem \eqref{1} as follows.
\begin{definition}\label{def1} Assume that $u_0, f\in L^\infty(\mathbb{R}^n)\cap L^1(\mathbb{R}^n)$. We say that $u = u(t, x)$ is a weak solution of problem \eqref{1}, if the following holds:
\begin{itemize}
\item[(a)] $u\in L_{loc}^\alpha((0,T)\times\mathbb{R}^n)\cap L_{loc}^p((0,T); W^{1,p}(\mathbb{R}^n)),$ $\nabla_x u\in L^\beta ((0,T)\times\mathbb{R}^n),$
    where $W^{1,p}(\mathbb{R}^n)$ is a Sobolev space.

\item[(b)] For all $0\leq \phi(t,x)\in C^1((0,T); C^1_0(\mathbb{R}^n)), \phi(T,x)=0, x\in\mathbb{R}^n,$
\begin{align*}
&\int\limits_0^T\int\limits_{\mathbb{R}^n}\left[\lambda|u(t,x)|^\alpha+\mu|\nabla_x u(t,x)|^\beta+f(x)\right]\phi(t,x)dxdt+\int\limits_{\mathbb{R}^n}u_0(x)\phi(0,x)dx\\&=\int\limits_0^T\int\limits_{\mathbb{R}^n}\left(|\nabla_x u(t,x)|^{p-2}\nabla_x u(t,x)\nabla_x\phi(t,x)-u(t,x)\phi_t(t,x)\right)dxdt.
\end{align*}
\end{itemize}
If $T=+\infty,$ then $u$ is called a global in time weak solution.
\end{definition}
Below we present our first main result.
\begin{theorem}[First critical exponents]\label{th1} Suppose that $u_0, f\in L^\infty(\mathbb{R}^n)\cap L^1(\mathbb{R}^n)$ and let $\beta,\alpha>\max\{1,p-1\},$ $\lambda>0,$ $\mu > 0.$
\begin{itemize}
\item[(i)] Let $f(x)>0$ for $p\neq 2$ or let $\int\limits_{\mathbb{R}^n}f(x)dx>0$ for $p=2.$ If
    $$\alpha\leq\frac{(p-1)n}{n-p}\,\,\,\, \text{or}\,\,\,\, \beta\leq\frac{(p-1)n}{n-1},$$ then the problem \eqref{1} admits no global weak solution;
\item[(ii)] Let $n > p.$ If
    $$\alpha>\frac{(p-1)n}{n-p}\,\,\,\, \text{and}\,\,\,\, \beta>\frac{(p-1)n}{n-1},$$ then the problem \eqref{1} admits positive global (stationary) classical solutions for sufficiently small $u_0>0$ and for some positive $f.$
\end{itemize}
\end{theorem}
\begin{remark}
It is easy to see from Theorem \ref{th1} that the critical exponents for the problem \eqref{1} are discontinuous.
Indeed, for $\beta>\frac{(p-1)n}{n-1},$ the critical exponent is the same as for the problem without the gradient term, that is, $\alpha_{cr}=\frac{(p-1)n}{n-p}$. However, the critical exponent for $\beta\leq\frac{(p-1)n}{n-1}$ is $\alpha_{cr}=\infty.$
\end{remark}
\begin{remark}
If $p=2$, Theorem \ref{th1} corresponds to Theorem 3 from \cite{Jleli1}. Part (i) of Theorem \ref{th1} also contains the critical cases of $\alpha$ and $\beta$ that remains new even when $p=2$. Therefore, it complements the result of the critical case, which was left open in \cite{Jleli1}.
\end{remark}
Then, for $n > p$, $\alpha>\frac{(p-1)n}{n-p}$ and $\beta>\frac{(p-1)n}{n-1},$ we investigate the effect of the forcing term $f$ on the critical behavior of problem \eqref{1}.

Let us define the sets
$$\mathcal{I}_r^+=\left\{f\in C(\mathbb{R}^n):\,f\geq 0, f(x)\geq C|x|^{-r}\,\,\,\text{for sufficiently large}\,|x|\right\}$$
and
$$\mathcal{I}_r^-=\left\{f\in C(\mathbb{R}^n):\,f> 0, f(x)\leq C|x|^{-r}\,\,\,\text{for sufficiently large}\,|x|\right\},$$
where $r<n$ and $C>0$ is an arbitrary constant.

Below we present the results for the second critical exponent in the sense of Lee and Ni \cite{Lee}.
\begin{theorem}[Second critical exponents]\label{th2} Let $u_0\in L^\infty(\mathbb{R}^n)\cap L^1(\mathbb{R}^n),$ $\alpha>\frac{(p-1)n}{n-p}$ and $\beta>\frac{(p-1)n}{n-1}.$
\begin{itemize}
\item[(i)] If $f\in \mathcal{I}^+_r$ and
    $$0<r<\max\left\{\frac{p\alpha}{\alpha-p+1}, \frac{\beta}{\beta-p+1}\right\},$$ then the problem \eqref{1} does not admit global positive weak solutions;
\item[(ii)] If
    $$\max\left\{\frac{p\alpha}{\alpha-p+1}, \frac{\beta}{\beta-p+1}\right\}\leq r<n,$$ then the problem \eqref{1} admits positive global (stationary) classical solutions for sufficiently small $u_0$ and for some positive $f\in\mathcal{I}_r^-.$
\item[(iii)] Suppose that $f\in \mathcal{I}^+_r$ and
    $r\leq 0,$ then the problem \eqref{1} does not admit any global positive weak solutions.
\end{itemize}
\end{theorem}
\begin{remark}
The results of Theorem \ref{th2} coincide with the results of Lee and Ni from \cite{Lee} for $\mu= 0$ and $p=2$.
Theorem \ref{th2} yields new results for all $p>\frac{2n}{n+1}$, including the case of the heat equation. It should be noted that the investigation of the second critical exponents of \eqref{1} in the case of $f=0$ remains open for all $p>\frac{2n}{n+1}$.
\end{remark}

\subsection{Historical background}\label{Hist} When $p=2,$ $\lambda=1,$ $\mu=0$ and $f=0,$ Fujita (see \cite{Fujita}) showed that the critical  exponent of problem \eqref{1} is $$\alpha_{cr}=1+\frac{2}{n},$$ which plays an important role in studying the existence and nonexistence of global solutions. That is:
\begin{itemize}
\item[(i)] if $1 <\alpha< 1+\frac{2}{n},$ then for positive initial values, problem \eqref{1} admits no global positive solutions;
\item[(ii)] if $\alpha > 1+\frac{2}{n},$ then for sufficiently small positive initial values, problem \eqref{1} admits positive global solutions.
\end{itemize}

In \cite{Galak}, Galaktionov extended Fujita's results for $\lambda=1,$ $\mu=0$ and $f=0$. Then the critical exponent of the problem \eqref{1} has the form
$$\alpha_{cr}=p-1+\frac{p}{n}.$$

The problem \eqref{1} with $\lambda=1,$ $\mu=0$ and $p=2$ was considered by Bandle-Levine-Zhang in \cite{Bandle}, and it was proven that:
\begin{itemize}
\item[(i)] if $1<\alpha<\alpha_{cr},\,u_0\geq 0$ and $\int_{\mathbb{R}^n}f(x) dx>0$, where \begin{equation*}\alpha_{cr}=\left\{\begin{array}{l}
\infty,\,\,\,\text{if}\,\,\,n=1,2, \\
\displaystyle\frac{n}{n-2},\,\,\,\text{if}\,\,\,n\geq3,\end{array}\right.\end{equation*} then the solution of problem \eqref{1} blows up in finite time;
\item[(ii)] if $n\geq 3$, $\alpha=\alpha_{cr}$, $\large\displaystyle\int_{\mathbb{R}^n}f(x)dx >0,$ and $$f(x)={O}(|x|^{-\varepsilon-n})\,\,\,\, \text{as} \,\,\,|x|\to\infty$$ for some $\varepsilon>0$, and either $u\geq0$ or
$$\int_{|x|>R}\frac{f^-(y)}{|x-y|^{n-2}}dy=\frac{o(1)}{|x|^{n-2}},\, f^-=\max\{-f, 0\},$$ when $R$ is enough large, then problem \eqref{1} has no global solutions;
\item[(iii)] if $n\geq 3$ and $\alpha> \alpha_{cr}$, then there exist $u_0\geq0$ and $f>0$ such that the problem \eqref{1} admits global solutions.
\end{itemize}
It should be noted that in \cite{Torebek}, the case (ii) has been improved without conditions on the asymptotic behavior of the function $f(x)$.

The results of \cite{Bandle} were extended in \cite{Zeng} for $\lambda=1,$ $\mu=0$ and $p>\frac{2n}{n+1},$ and it was demonstrated that the critical exponent takes the form
\begin{equation*}\alpha_{cr}=\left\{\begin{array}{l}
\infty,\,\,\,\text{if}\,\,\,n\leq p, \\
\displaystyle\frac{(p-1)n}{n-p},\,\,\,\text{if}\,\,\,n\geq p.\end{array}\right.\end{equation*}

The equation \eqref{1} coincides with the viscous Hamilton-Jacobi equation when $\lambda=0$.
The large-time behavior of positive solutions to viscous Hamilton-Jacobi equations with and without forcing terms was studied by many authors (see for example \cite{Amal, Weissler, Souplet1, Souplet2, Souplet3} and references therein). In particular, when $p=2,$ $\mu>0$ and $f=0$ it is known that (see \cite{Souplet1}):
\begin{itemize}
\item[(i)] if $1 < \beta \leq\frac{n+2}{n+1},$ then $I_\infty:=\lim\limits_{t\rightarrow+\infty}\|u(t,\cdot)\|_{L^1(\Omega)} = \infty$ for all $u_0>0$;
\item[(ii)] if $\frac{n+2}{n+1} <\beta< 2,$ then both $I_\infty =\infty$ and $I_\infty < \infty$ occur;
\item[(iii)] if $\beta \geq 2,$ then $I_\infty < \infty$ for all $u_0>0$.
\end{itemize}

In \cite{Jleli1}, Jleli-Samet-Souplet considered problem \eqref{1} for $p=2$ and obtained the following results:

(a) In the case $f(x)=0$
\begin{itemize}
\item[(1a)] if $1<\alpha\leq 1+\frac{2}{n}$ or $1<\beta\leq 1+\frac{1}{n+1},$ then for any $u_0>0$ the solution of problem \eqref{1}  blows up in finite time;
\item[(2a)] if $\alpha> 1+\frac{2}{n}$ and $\beta> 1+\frac{1}{n+1},$ then the problem \eqref{1} admits positive global solutions for suitably small initial data.
\end{itemize}

(b) In the case $f(x)\neq0$
\begin{itemize}
\item[(1b)] if $1<\alpha< \frac{n}{n-2}$ or $1<\beta< \frac{n}{n-1},$ then for any $u_0\in BC^1(\mathbb{R}^n)$ and $f\in BC^1(\mathbb{R}^n)\cap L^1(\mathbb{R}^n)$ such that $\int_{\mathbb{R}^n}f(x) dx>0,$ the solution of problem \eqref{1} blows up in finite time;
\item[(2b)] if $\alpha> \frac{n}{n-2}$ and $\beta> \frac{n}{n-1},$ then the problem \eqref{1} admits positive global solutions for suitably small initial data.
\end{itemize}

The presence of a gradient term causes the phenomenon of critical exponent discontinuity, as shown in the results above. Namely, as long as $\beta> \frac{n}{n-1},$ the critical exponent remains the same as for the problem without the gradient term, i.e., $\alpha_{cr} =\frac{n}{n-2}$ and when $\beta<\frac{n}{n-1}$, the critical exponent is $\alpha_{cr} = \infty.$

It should be noted that when $f(x)\neq 0$, the critical cases $$\alpha=\frac{n}{n-2}\,\,\,\,\text{and}\,\,\,\,\beta=\frac{n}{n-1},$$ were not investigated and were left \textbf{open} in \cite{Jleli1}.

It is also worth noting that the existence of global and blow-up solutions to problem \eqref{1} for $p=2$ were previously studied in \cite{Weissler1, Zaag}.

Recently, in \cite{Lu}, Lu-Zhang extended the result of Jleli-Samet-Souplet to the case $p$-Laplace heat equation (i.e. problem \eqref{1} for $f=0$) and discovered a critical exponent of the form
$$\alpha_{cr}=p-1+\frac{p}{n}\,\,\,\,\text{and}\,\,\,\,\beta_{cr}=p-1+\frac{1}{n+1}.$$ Also, the discontinuity of the critical exponents is proved.

Motivated by the preceding results, we extend the results of Jleli-Samet-Souplet for degenerate equation \eqref{1} with $f(x)\neq 0$. In addition, we study the critical cases, and in particular (when $p=2$), we will fill the gap the results of \cite{Jleli1}.

\section{Proof of Theorem \ref{th1}}
In this section, we present the proofs of the main results step by step for each case.
\subsection{The case $\alpha<\frac{(p-1)n}{n-p}$ or $\beta<\frac{(p-1)n}{n-1}$}\label{S1}
The case $$\alpha<\frac{(p-1)n}{n-p}$$ is obvious, because the positive solution $u_\lambda$ to problem \eqref{1} without a gradient term blows up in finite time (see \cite{Zeng} and \cite{Bandle}). As $u \geq u_\lambda,$ according to the comparison principle, the solution to problem \eqref{1} also blows up when $\alpha<\frac{(p-1)n}{n-p}.$

We now focus on the case $\beta<\frac{(p-1)n}{n-1}.$ Assume that there is a global weak solution to problem \eqref{1} for $\beta<\frac{(p-1)n}{n-1}.$
According to Definition \ref{def1}, for a weak solution to problem \eqref{1} we have
\begin{align*}
&\int\limits_0^T\int\limits_{\mathbb{R}^n}\left[\lambda|u(t,x)|^\alpha+\mu|\nabla_x u(t,x)|^\beta+f(x)\right]\phi(t,x)dxdt+\int\limits_{\mathbb{R}^n}u_0(x)\phi(0,x)dx\\&=\int\limits_0^T\int\limits_{\mathbb{R}^n}\left(|\nabla_x u(t,x)|^{p-2}\nabla_x u(t,x)\nabla_x\phi(t,x)-u(t,x)\phi_t(t,x)\right)dxdt.
\end{align*}
Due to $\varepsilon$-Young inequality we have
\begin{align*}
\int\limits_0^T\int\limits_{\mathbb{R}^n}u\phi_tdxdt&=\int\limits_0^T\int\limits_{\mathbb{R}^n}u\phi^{\frac{1}{\alpha}}\phi_t\phi^{-\frac{1}{\alpha}}dxdt \\&\leq \varepsilon_1\int\limits_0^T\int\limits_{\mathbb{R}^n} |u|^\alpha\phi dxdt +C_1(\varepsilon_1)\int\limits_0^T\int\limits_{\mathbb{R}^n}\frac{|\phi_t|^{\frac{\alpha}{\alpha-1}}}{|\phi|^\frac{1}{\alpha-1}}dxdt,
\end{align*}
and
\begin{align*}
\int\limits_0^T\int\limits_{\mathbb{R}^n}|\nabla_x u|^{p-2}\nabla_x u\nabla_x\phi dxdt&=\int\limits_0^T\int\limits_{\mathbb{R}^n}|\nabla_x u|^{p-2}\nabla_x u\phi^{\frac{p-1}{\beta}}\nabla_x\phi\phi^{-\frac{p-1}{\beta}} dxdt \\&\leq \int\limits_0^T\int\limits_{\mathbb{R}^n}|\nabla_x u|^{p-1}\phi^{\frac{p-1}{\beta}}|\nabla_x\phi|\phi^{-\frac{p-1}{\beta}} dxdt \\&\leq \varepsilon_2\int\limits_0^T\int\limits_{\mathbb{R}^n}|\nabla_x u|^\beta\phi dxdt+C_2(\varepsilon_2)\int\limits_0^T\int\limits_{\mathbb{R}^n}\frac{|\nabla_x \phi|^{\frac{\beta}{\beta-p+1}}}{|\phi|^\frac{p-1}{\beta-p+1}}dxdt,
\end{align*} where $C_1$ and $C_2$ are some positive constants.

Then choosing $\varepsilon_1=\frac{\lambda}{2}$ and $\varepsilon_2=\frac{\mu}{2}$ one can obtain
\begin{equation}\label{4}
\begin{split}
\frac{1}{2}\int\limits_0^T\int\limits_{\mathbb{R}^n}&\left(\lambda|u|^\alpha+\mu|\nabla_x u|^\beta\right)\phi dxdt+\int\limits_0^T\int\limits_{\mathbb{R}^n}f(x)\phi dxdt+\int\limits_{\mathbb{R}^n}u_0(x)\phi(0,x)dx\\&\leq C_1\int\limits_0^T\int\limits_{\mathbb{R}^n}\frac{|\phi_t|^{\frac{\alpha}{\alpha-1}}}{|\phi|^\frac{1}{\alpha-1}}dxdt +C_2\int\limits_0^T\int\limits_{\mathbb{R}^n}\frac{|\nabla_x \phi|^{\frac{\beta}{\beta-p+1}}}{|\phi|^\frac{p-1}{\beta-p+1}}dxdt. \end{split}
\end{equation}
Let $\phi$ satisfy
\begin{equation}\label{5}
\phi(t,x)=\Psi^{\frac{\alpha}{\alpha-1}}\left(\frac{t}{T}\right)\Phi^{\frac{\beta}{\beta-p+1}}\left(\frac{|x|}{T^\kappa}\right), \, \Psi\in C^1_0(\mathbb{R}_+),\, \Phi\in C_0^1(\mathbb{R}),
\end{equation}
and
$$\Phi(s) =
\begin{cases}
1, & \text{if $ 0\leq s \leq 1/2 $,} \\
\searrow, & \text{if $1/2< s\leq 1$,}\\
0, & \text{if $s>  1$,}
\end{cases}\,\,\,
\Psi(s) =
\begin{cases}
1, & \text{if $ 0\leq s \leq 1/2 $,} \\
\searrow, & \text{if $1/2< s\leq 1$,}\\
0, & \text{if $s>  1$.}
\end{cases}$$
By the properties of $\Phi$ and $\Psi,$ and changing $t=T\tau,\, x=T^\kappa y$ we arrive at
\begin{align*}\int\limits_0^T\int\limits_{\mathbb{R}^n}{|\phi_t|^{\frac{\alpha}{\alpha-1}}}{|\phi|^{-\frac{1}{\alpha-1}}}dxdt&\leq
CT^{\kappa n+1-\frac{\alpha}{\alpha-1}}\int\limits_{\frac{1}{2}}^1|\Psi'(\tau)|^{{\frac{\alpha}{\alpha-1}}}d\tau\int\limits_{|y|\leq 1}\Phi^{\frac{\beta}{\beta-p+1}}(|y|)dy,
\end{align*}
and
\begin{align*}
\int\limits_0^T\int\limits_{\mathbb{R}^n}\frac{|\nabla_x \phi|^{\frac{\beta}{\beta-p+1}}}{|\phi|^\frac{p-1}{\beta-p+1}}dxdt&\leq
CT^{\kappa n+1-\kappa\frac{\beta}{\beta-p+1}}\int\limits_0^1\Psi^{\frac{\alpha}{\alpha-1}}(\tau)d\tau\int\limits_{\frac{1}{2}\leq|y|\leq 1}|\Phi'(|y|)|^{\frac{\beta}{\beta-p+1}}dy.
\end{align*}
Then choosing $\kappa=\frac{\alpha(\beta-p+1)}{\beta(\alpha-1)}$ we deduce that
\begin{equation*}
\begin{split}
\int\limits_0^T\int\limits_{\mathbb{R}^n}f(x)\phi dxdt+\int\limits_{\mathbb{R}^n}u_0(x)\phi(0,x)dx\leq CT^{\frac{\alpha(\beta-p+1)}{\beta(\alpha-1)} n+1-\frac{\alpha}{\alpha-1}}.\end{split}
\end{equation*}
It is obvious that
\begin{equation*}
\begin{split}
\int\limits_0^T\int\limits_{\mathbb{R}^n}f(x)\phi dxdt=T\int\limits_0^1\Psi^{\frac{\alpha}{\alpha-1}}(\tau)d\tau\int\limits_{\mathbb{R}^n}f(x)\Phi^{\frac{\beta}{\beta-p+1}}\left(\frac{|x|}{T^\kappa}\right) dx.\end{split}\end{equation*}
Consequently
\begin{equation}\label{6}
\begin{split}
\int\limits_{\mathbb{R}^n}f(x)\Phi^{\frac{\beta}{\beta-p+1}}\left(\frac{|x|}{T^\kappa}\right) dx\leq C_3T^{-1}\int\limits_{\mathbb{R}^n}|u_0(x)|\phi(0,x)dx+C_4T^{\frac{\alpha(\beta-p+1)}{\beta(\alpha-1)} n-\frac{\alpha}{\alpha-1}},\end{split}
\end{equation} where $C_3$ and $C_4$ are some positive constants.

Let $p\neq 2$ and $f>0,$ then there exist $\delta>0$ such that $$\int\limits_{\mathbb{R}^n}f(x)\Phi^{\frac{\beta}{\beta-p+1}}\left(\frac{|x|}{T^\kappa}\right) dx\geq \int\limits_{B_{\frac{1}{2}}(0)}f(x)dx>\delta,$$ hence from \eqref{6} we have that
\begin{equation}\label{6*}
\begin{split}
\delta\leq C_3T^{-1}\int\limits_{\mathbb{R}^n}|u_0(x)|\phi(0,x)dx+C_4T^{\frac{\alpha(\beta-p+1)}{\beta(\alpha-1)} n-\frac{\alpha}{\alpha-1}},\end{split}\end{equation}

As $\beta<\frac{(p-1)n}{n-1}$ it follows that $\frac{\alpha}{\alpha-1}\left(\frac{(\beta-p+1)}{\beta} n-1\right)<0.$ Hence, passing to the limit as $T \rightarrow +\infty$ in \eqref{6*}, one obtains
$\delta\leq 0,$ which contradicts the fact that $\delta>0.$

Let $p=2$ and $\int\limits_{\mathbb{R}^n}f(x) dx>0.$ As $\beta<\frac{n}{n-1}$ it follows that $\frac{\alpha}{\alpha-1}\left(\frac{(\beta-1)}{\beta} n-1\right)<0.$ Hence, passing to the limit as $T \rightarrow +\infty$ in \eqref{6}, one obtains
\begin{equation*}
\lim\limits_{T \rightarrow +\infty}\int\limits_{\mathbb{R}^n}f(x)\Phi^{\frac{\beta}{\beta-p+1}}\left(\frac{|x|}{T^\kappa}\right) dx=\int\limits_{\mathbb{R}^n}f(x) dx\leq 0,
\end{equation*} which contradicts the fact that $\int\limits_{\mathbb{R}^n}f(x) dx>0.$

\subsection{The case $\alpha=\frac{(p-1)n}{n-p}$ or $\beta=\frac{(p-1)n}{n-1}$}
The result for $$\alpha=\frac{(p-1)n}{n-p}$$ follows immediately from \cite{Zeng, Bandle}.

Let $\beta=\frac{(p-1)n}{n-1}.$ Assume that for $\beta=\frac{(p-1)n}{n-1},$ the global solution to problem \eqref{1} exists.

Repeating the procedure as in subsection \ref{S1} from \eqref{4}, we have
\begin{align*}
\int\limits_0^T\int\limits_{\mathbb{R}^n}f(x)\phi dxdt&+\int\limits_{\mathbb{R}^n}u_0(x)\phi(0,x)dx\\&\leq C_1\int\limits_0^T\int\limits_{\mathbb{R}^n}\frac{|\phi_t|^{\frac{\alpha}{\alpha-1}}}{|\phi|^\frac{1}{\alpha-1}}dxdt +C_2\int\limits_0^T\int\limits_{\mathbb{R}^n}\frac{|\nabla_x \phi|^{\frac{\beta}{\beta-p+1}}}{|\phi|^\frac{p-1}{\beta-p+1}}dxdt,
\end{align*} where $C_1$ and $C_2$ are some positive constants.

Let $\phi$ satisfy
\begin{equation*}
\phi(t,x)=\Psi^{\frac{\alpha}{\alpha-1}}\left(\frac{t}{T}\right) \Phi^{\frac{\beta}{\beta-p+1}}\left(\frac{\ln\left(\frac{|x|}{{T^\theta}}\right)}{\ln\left({T^\theta}\right)}\right), \, T\gg1,\, \theta>0,
\end{equation*}
where $\Psi\in C^1_0(\mathbb{R}_+)$ and $\Phi\in C_0^1(\mathbb{R})$ satisfy
$$\Psi(s) =
\begin{cases}
1, & \text{if $ 0\leq s \leq 1/2 $,} \\
\searrow, & \text{if $1/2< s< 1$,}\\
0, & \text{if $s\geq  1$,}
\end{cases}$$
and
\begin{equation*}
\Phi(s)=\left\{\begin{array}{l}
1,\,\,\,\,\,\text{if}\,\,-\infty< s\leq0,\\
\searrow, \,\,\text{if} \,\,\,{0< s< 1,}\\
0,\,\,\,\,\,\text{if}\,\,s\geq1.\end{array}\right.
\end{equation*}

Using the properties of $\Phi$ and $\Psi,$ as well as changing $t=T\tau,\, x={T^\theta} y,$ and applying
$$\left|\nabla \Phi\left(\frac{\ln\left(\frac{|x|}{{T^\theta}}\right)}{\ln\left({T^\theta}\right)}\right)\right|\leq \frac{C}{|x|\ln T},\,\,C>0,$$
we have
\begin{align*}\int\limits_0^T\int\limits_{\mathbb{R}^n}\frac{|\phi_t|^{\frac{\alpha}{\alpha-1}}}{|\phi|^\frac{1}{\alpha-1}}dxdt&\leq
CT^{2\theta n+1- \frac{\alpha}{\alpha-1}}
\end{align*}
and
\begin{align*}
\int\limits_0^T\int\limits_{\mathbb{R}^n}\frac{|\nabla_x \phi|^{n}}{|\phi|^{n-1}}dxdt&\leq
CT(\ln T)^{1-n}.
\end{align*}
Hence
\begin{equation*}
\begin{split}
\int\limits_{\mathbb{R}^n}f(x)\Phi^{n}\left(\frac{\ln\left(\frac{|x|}{{T^\theta}}\right)}{\ln\left({T^\theta}\right)}\right) dx&\leq C_1\left(T^{2\theta n- \frac{\alpha}{\alpha-1}}+(\ln T)^{1-n}\right)\\&+C_2T^{-1}\int\limits_{\mathbb{R}^n}u_0(x)\phi(0,x)dx,\end{split}
\end{equation*} where $C_1$ and $C_2$ are some positive constants.

As
$$\lim\limits_{T\rightarrow+\infty}\Phi\left(\frac{\ln\left(\frac{|x|}{{T^\theta}}\right)}{\ln\left({T^\theta}\right)}\right)=\Phi(-1)=1,$$
then choosing $\theta < \frac{\alpha}{2n(\alpha-1)},$ one can obtain
\begin{equation*}
\begin{split}
\lim\limits_{T\rightarrow+\infty}\int\limits_{\mathbb{R}^n}f(x) \Phi^{n}\left(\frac{\ln\left(\frac{|x|}{{T^\theta}}\right)}{\ln\left({T^\theta}\right)}\right) dx=\int\limits_{\mathbb{R}^n}f(x)dx\leq 0,\end{split}
\end{equation*} which contradicts the fact that $\int\limits_{\mathbb{R}^n}f(x) dx>0.$

\subsection{The case $\alpha>\frac{(p-1)n}{n-p}$ and $\beta>\frac{(p-1)n}{n-1}$} Let
\begin{equation*}
v(x)=\epsilon\left(1+|x|^{\frac{p}{p-1}}\right)^{-m},
\end{equation*} where $\epsilon>0$ and $$\max\left\{\frac{p-1}{\alpha-p+1},\frac{(p-1)(p-\beta)}{p(\beta-p+1)}\right\}<m<\frac{n-p}{p}.$$

Let us compute
\begin{align*}|\nabla v(x)|&=m\epsilon\frac{p}{p-1}|x|^{\frac{1}{p-1}}\left(1+|x|^{\frac{p}{p-1}}\right)^{-m-1}\\&\leq m\epsilon\frac{p}{p-1}\left(1+|x|^{\frac{p}{p-1}}\right)^{-m-1+\frac{1}{p}},\end{align*}
and
\begin{align*}-\Delta_p v(x)&=n\left(\frac{\epsilon mp}{p-1}\right)^{p-1}\left(1+|x|^{\frac{p}{p-1}}\right)^{-(m+1)(p-1)}\\&- (m+1)p\left(\frac{\epsilon mp}{p-1}\right)^{p-1}|x|^{\frac{p}{p-1}}\left(1+|x|^{\frac{p}{p-1}}\right)^{-(m+1)(p-1)-1}\\& >\left(\frac{\epsilon mp}{p-1}\right)^{p-1}(n-p-mp)\left(1+|x|^{\frac{p}{p-1}}\right)^{-(m+1)(p-1)}.\end{align*}
Then
\begin{align*}&-\Delta_p v-v^\alpha-|\nabla v|^\beta\\&>\left(\frac{\epsilon mp}{p-1}\right)^{p-1}(n-p-mp)\left(1+|x|^{\frac{p}{p-1}}\right)^{-(m+1)(p-1)}\\&-\epsilon^\alpha\left(1+|x|^{\frac{p}{p-1}}\right)^{-m\alpha} - \left(\frac{\epsilon mp}{p-1}\right)^\beta\left(1+|x|^{\frac{p}{p-1}}\right)^{-(m+1)\beta+\frac{\beta}{p}}.
\end{align*}
As $m>\max\left\{\frac{p-1}{\alpha-p+1},\frac{(p-1)(p-\beta)}{p(\beta-p+1)}\right\},$ $\alpha>\frac{(p-1)n}{n-p}$ and $\beta>\frac{(p-1)n}{n-1},$ it follows that
$$\left(1+|x|^{\frac{p}{p-1}}\right)^{-m\alpha}<\left(1+|x|^{\frac{p}{p-1}}\right)^{-(m+1)(p-1)}$$
and
$$\left(1+|x|^{\frac{p}{p-1}}\right)^{-(m+1)\beta+\frac{\beta}{p}}<\left(1+|x|^{\frac{p}{p-1}}\right)^{-(m+1)(p-1)}.$$
Hence
\begin{align*}-\Delta_p v-v^\alpha-|\nabla v|^\beta>\left(\frac{ \epsilon mp}{p-1}\right)^{p-1}M_\epsilon\left(1+|x|^{\frac{p}{p-1}}\right)^{-(m+1)(p-1)},
\end{align*} where $M_\epsilon=(n-p-mp)-\epsilon^{\alpha-p+1}\left(\frac{p-1}{mp}\right)^{p-1} - \left(\frac{\epsilon mp}{p-1}\right)^{\beta-p+1}.$

Then, choosing $\epsilon$ sufficiently small, we obtain
\begin{align*}-\Delta_p v-v^\alpha-|\nabla v|^\beta>0,\,x\in\mathbb{R}^n.\end{align*}
Setting $f(x):= -\Delta_p v-v^\alpha-|\nabla v|^\beta>0,$ we can deduce that the function $v$ is a stationary
positive supersolution to problem \eqref{1}, i.e $u(t,x)\leq v(x)$ for all $(t,x)\in\mathbb{R}_+\times\mathbb{R}^n.$

\section{Proof of Theorem \ref{th2}}
In this section, we present the proofs step by step for each case.
The proof is based on the approach of the previous section, but slightly modified in accordance with the assumptions of Theorem \ref{th2}.
\subsection{The case $r<\max\left\{\frac{p\alpha}{\alpha-p+1}, \frac{\beta}{\beta-p+1}\right\}$}\label{S2}
Since some of the a priori estimates used in the proof of Theorem \ref{th1} still work, it remains for us to consider some additional factors.

Let $r<\frac{p\alpha}{\alpha-p+1}.$ Then, repeating the technique from \cite{Zeng} with the combination of calculations from subsection \ref{S1}, we come to the conclusion that for $$r<\frac{p\alpha}{\alpha-p+1},$$ there is no global in time nontrivial solutions of problem \eqref{1}.

Let $r<\frac{\beta}{\beta-p+1}.$
Using the same proof scheme as in subsection \ref{S1}, from \eqref{6} we have
\begin{equation}\label{7}
\begin{split}
\int\limits_{\mathbb{R}^n}f(x)\Phi^{\frac{\beta}{\beta-p+1}}\left(\frac{|x|}{T^\kappa}\right) dx&\leq C_1T^{-1}\int\limits_{\mathbb{R}^n}|u_0(x)|\phi(0,x)dx\\& +C_2(T^{\kappa n-\frac{\alpha}{\alpha-1}}+T^{\kappa n-\frac{\kappa\beta}{\beta-p+1}}).\end{split}
\end{equation}
Suppose that $\Phi$ has the following property:
$$\Phi(s) =
\begin{cases}
1, & \text{if $ 0\leq s < 1 $,} \\
\searrow, & \text{if $1\leq s< 2$,}\\
0, & \text{if $s\geq  2$.}
\end{cases}$$
Let us estimate the left side of \eqref{7}. As $f\in\mathcal{I}^+_r,$ we have
\begin{align*}\int\limits_{\mathbb{R}^n}f(x)\Phi^{\frac{\beta}{\beta-p+1}}\left(\frac{|x|}{T^\kappa}\right) dx&= \int\limits_{|x|< 2T^\kappa}f(x)\Phi^{\frac{\beta}{\beta-p+1}}\left(\frac{|x|}{T^\kappa}\right) dx\\&\geq \int\limits_{|x|< T^\kappa}f(x)\Phi^{\frac{\beta}{\beta-p+1}}\left(\frac{|x|}{T^\kappa}\right) dx\\&=\int\limits_{|x|< T^\kappa}f(x) dx\geq \int\limits_{\frac{T^\kappa}{2}<|x|<T^\kappa}f(x) dx\\&\geq C\int\limits_{\frac{T^\kappa}{2}<|x|<T^\kappa}|x|^{-r} dx=CT^{(n-r)\kappa},\,\,C>0.
\end{align*}
Then, combining this with \eqref{7}, we obtain
\begin{equation}\label{8}
CT^{(n-r)\kappa}\leq C_1T^{-1}\int\limits_{\mathbb{R}^n}|u_0(x)|dx+C_2(T^{\kappa n-\frac{\alpha}{\alpha-1}}+T^{\kappa n-\frac{\kappa\beta}{\beta-p+1}}).\end{equation}

Choosing $\kappa=\frac{\alpha(\beta-p+1)}{\beta(\alpha-1)}$ in \eqref{8} we deduce that
\begin{equation}\label{9}
C\leq C_1T^{-1-(n-r)\kappa}\int\limits_{\mathbb{R}^n}|u_0(x)|dx+C_2T^{\frac{\alpha}{\alpha-1}{\left(\frac{\beta-p+1}{\beta} r-1\right)}}.\end{equation}
As $r<\frac{\beta}{\beta-p+1}$ it follows that $(\beta-p+1)r-\beta<0.$ Hence, passing to the limit as $T \rightarrow +\infty$ in \eqref{9}, one obtains
\begin{equation*}
C\leq 0,
\end{equation*} which contradicts the fact that $C>0.$
\subsection{The case $\max\left\{\frac{p\alpha}{\alpha-p+1}, \frac{\beta}{\beta-p+1}\right\}\leq r<n$} We consider the function
\begin{equation*}
v(x)=\epsilon\left(1+|x|^{\frac{p}{p-1}}\right)^{-m},
\end{equation*} where $\epsilon>0$ and $$\frac{r-p}{p}<m<\frac{n-p}{p}.$$
The set of values of $m$ satisfying above inequality is nonempty, as $r<n$. Furthermore,  $m$ also satisfies $\max\left\{\frac{p-1}{\alpha-p+1},\frac{(p-1)(p-\beta)}{p(\beta-p+1)}\right\}<m<\frac{n-p}{p}.$ Consequently, from the proof of part (ii) of Theorem \ref{th1}, it follows that
$$f(x)= -\Delta_p v-v^\alpha-|\nabla v|^\beta>0,\,x\in\mathbb{R}^n,$$ for sufficiently small $\epsilon>0.$ Hence $v(x)$ is a positive stationary (super)solution of \eqref{1}.

On the other hand, using elementary calculations, it is easy to show
\begin{align*}f(x)&\leq n\left(\frac{\epsilon mp}{p-1}\right)^{p-1}\left(1+|x|^{\frac{p}{p-1}}\right)^{-(m+1)(p-1)} \\&\leq C\left(1+|x|^{\frac{p}{p-1}}\right)^{-\frac{r(p-1)}{p}}\\& \leq C|x|^{-r},\,\,C>0,\end{align*}
which asserts $f\in\mathcal{I}^-_r.$ The proof is complete.

\subsection{The case $r\leq 0$} The proof follows from the results of Subsection \ref{S1}. Indeed, let $r\leq 0$, then from \eqref{8} we obtain the following
\begin{equation*}
C\leq C_1T^{-1-(n-r)\kappa}+C_2(T^{\kappa r-\frac{\alpha}{\alpha-1}}+T^{\kappa r-\frac{\kappa\beta}{\beta-p+1}}).\end{equation*}
Since $r\leq 0$ it follows that $$-1-(n-r)\kappa<0,$$ $$\kappa r-\frac{\alpha}{\alpha-1}<0,$$ and $$\kappa r-\frac{\kappa\beta}{\beta-p+1}<0.$$ Consequently, passing to the limit as $T \rightarrow +\infty$ in last inequality, one obtains
\begin{equation*}
0< C\leq 0,
\end{equation*} which contradicts the fact that $C>0.$

\section*{Declaration of competing interest}
	The Author declare that there is no conflict of interest

\section*{Data Availability Statements} The manuscript has no associated data

\section*{Acknowledgments}
The author would like to thank the reviewers for their valuable comments
and remarks.

\end{document}